# Machine Learning and Deep Learning in Computational Finance: A Systematic Review

**Soufiane El Amine El Alami**
Euromed University of Fes, UEMF, Morocco
Email: soufiane.el-amine-el-alami@insa.ueuromed.org
ORCID iD: https://orcid.org/0009-0005-1770-4076
*Corresponding Author

**Abderazzak Mouiha**
Euromed University of Fes, UEMF, Morocco
Email: a.mouiha@ueuromed.org
ORCID iD: https://orcid.org/0000-0002-7034-7267

**Abdelatif Hafid**
FSTH, Abdelmalek Essaadi University, Morocco
Email: a.hafid@uae.ac.ma
ORCID iD: https://orcid.org/0000-0002-1377-1387

**Ahmed El Hilali Alaoui**
Euromed University of Fes, UEMF, Morocco
Email: a.elhilali-alaoui@ueuromed.org

**Abstract:** This systematic review examines how machine learning (ML) and deep learning (DL) have transformed forecasting, decision-making, and financial modelling, promoting innovation and efficiency in financial systems. Following PRISMA 2020 guidelines, we analyze 22 peer-reviewed and open-access articles (2024–2026) indexed in Scopus, applying ML and DL models across credit risk prediction, cryptocurrency, asset pricing, and macroeconomic policy modeling. The most used models include Random Forest, XG-Boost, Support Vector Machine, Long Short-Term Memory (LSTM), Bidirectional LSTM, Convolutional Neural Network (CNN), and hybrid or ensemble approaches combining statistical and AI methods. ML and DL techniques outperform traditional models by capturing nonlinear dependencies and enhancing predictive accuracy, while explainable AI methods (e.g., SHAP and feature importance analysis) improve transparency and interpretability. Emerging trends include cross-domain applications and the integration of responsible AI in finance. Despite notable progress, challenges remain in interpretability, generalizability, and data quality. Overall, this review provides a comprehensive overview of AI-driven computational finance and outlines future research directions.

## 1. Introduction

In recent years, the combination of computational finance and artificial intelligence has become one of the most dynamic areas in economics. Machine learning (ML) and Deep learning (DL) have changed how financial modeling, forecasting, and decision-making are performed [1]. These methods can learn complex, nonlinear relationships in large datasets, allowing researchers and practitioners to make more accurate predictions and to improve financial decision processes. The growing availability of high-frequency and multidimensional financial data, together with advances in computing power, have made these techniques more efficient and accessible than ever before [2].

Computational finance has particularly benefited from the flexibility of ML and DL, which allow for modeling patterns that traditional statistical models cannot easily capture. Classical approaches like support vector machines, random forests, or XGBoost are often used for credit risk and asset pricing, while deep architectures such as LSTM or CNN are applied to time series forecasting and anomaly detection. Recent research also explores hybrid and ensemble methods that combine different algorithms to enhance prediction accuracy [3], [4]. For instance, hybrid ML models for volatility prediction have improved the management of financial risks, and deep unsupervised models have shown strong performance in anomaly detection in high-frequency trading environments [5].

The most recent literature demonstrates that ML and DL are being used across almost every area of finance [1]. Studies published in 2024 and 2025 show applications in credit scoring, stock market forecasting, cryptocurrency price prediction, risk management, and portfolio optimization [6], [7], [8], [9], [10]. Some works focus on specific financial systems such as the Turkish banking sector [11], the Australian equity market [12], or the Saudi stock market [13]. Others apply ML to macroeconomic or policy-related forecasting, insurance and actuarial models, or financial technology (FinTech) innovation. These diverse applications confirm that ML and DL are no longer limited to academic experiments but have become essential tools for real-world decision-making.

Despite this rapid growth, the literature remains fragmented. Many papers focus on a single asset class or a narrow use case, such as cryptocurrency prediction or credit card default risk [7], [14], while fewer provide a broad and systematic understanding of how ML and DL are used across financial domains [1]. There is also a lack of bibliometric synthesis connecting recent empirical evidence to methodological progress. For this reason, the present study conducts a bibliometric review of current contributions to computational finance that explicitly employ ML and DL methods.

The objectives of this review are threefold. First, it aims to identify and classify the most frequently used ML and DL models in computational finance. Second, it evaluates their empirical effectiveness compared to traditional statistical models (e.g., linear regression, logistic regression, ARIMA, VAR, and GARCH). Third, it highlights the methodological challenges and practical implications that researchers and professionals encounter when applying these models to real financial data. In addition, this review considers how these approaches contribute to emerging areas such as sustainable finance, FinTech entrepreneurship, and data-driven policy forecasting.

Guided by these objectives, the review addresses the following research questions:

RQ1: What are the most commonly used ML and DL models in computational finance?
RQ2: How do these models perform compared to traditional statistical methods?
RQ3: What are the main challenges and limitations in applying ML and DL to financial forecasting?

By answering these questions, this work provides an updated and structured view of the field. It offers evidence on how recent advances in artificial intelligence continue to reshape modern finance, both in theory and in practice.

## 2. Methods (Following PRISMA Guidelines)

This section outlines the systematic review methodology following PRISMA 2020, describes the Scopus search strategy and the inclusion and exclusion criteria, and eventually explains the data extraction process.

### *2.1 Search Strategy*

This systematic review follows the PRISMA 2020 guidelines [15]. The objective is to identify recent peer-reviewed studies that apply ML and DL methods in the field of computational finance and financial forecasting.

The literature search was conducted using Scopus, which was selected as the main database for three main reasons. First, Scopus provides comprehensive coverage of peer-reviewed journals in economics, finance, and data science. Second, it ensures that only high-quality academic publications are included. Third, Scopus is the only full-access database available through my institutional credentials, which guarantees consistent and verifiable access to all selected papers.

The updated search query used in this review was refined to capture the most recent literature between 2024 and 2026. The exact query applied on ScienceDirect (Scopus-indexed journals) was: ("machine learning" OR "deep learning") AND "finance".

The search was further restricted using the following filters:

- Years: 2024–2026
- Article type: Research articles only
- Access type: Open access
- Subject area: Economics, Econometrics, and Finance

This combination of filters ensured that the analysis focused on recent, peer-reviewed, and openly accessible research within the economic and financial domains. A total of 480 records were initially retrieved (as shown in Figure 1).

The inclusion criteria were:

- Peer-reviewed journal articles written in English
- Studies explicitly applying ML or DL methods in finance

- Research with an empirical design or model evaluation
- Availability of full text under an open access license

The exclusion criteria were:

- Preprints, conference papers, or technical reports
- Studies outside the field of economics or finance
- Articles without quantitative or forecasting components
- Duplicates and incomplete metadata

*2.2 Study Selection (PRISMA Flow Chart)*

The study selection process followed the PRISMA framework and was carried out in four distinct stages:

1. Identification: 480 records were retrieved from the Scopus search.
2. Screening: Titles and abstracts were reviewed to remove irrelevant studies.
3. Eligibility: 72 full-text articles were examined in detail against the inclusion criteria.
4. Inclusion: Finally, 22 articles were selected for the bibliometric and content analysis.

Each selected study was checked manually to confirm that it contained clear references to ML or DL applications in finance, with reproducible methods and results. The selection process is summarized in the PRISMA flow chart (Figure 1).

For each study, relevant information was extracted, including:

- authorship and year of publication,
- journal and publisher,
- ML/DL models used,
- dataset characteristics (size, frequency, asset type),
- performance metrics (RMSE, MAE, Accuracy, F1-score, etc.),
- and key findings or limitations.

The final dataset of 22 studies represents a comprehensive sample of the most recent research applying machine and deep learning to financial forecasting, asset pricing, risk prediction, and decision-making.

Figure 1 shows the process used to identify, screen, and include studies in the systematic review. A total of 480 records were retrieved from Scopus after applying the search filters (2024–2026, research articles, open access, and the subject area of Economics, Econometrics, and Finance). After screening and eligibility assessment, 22 articles met the inclusion criteria and were retained for final analysis.

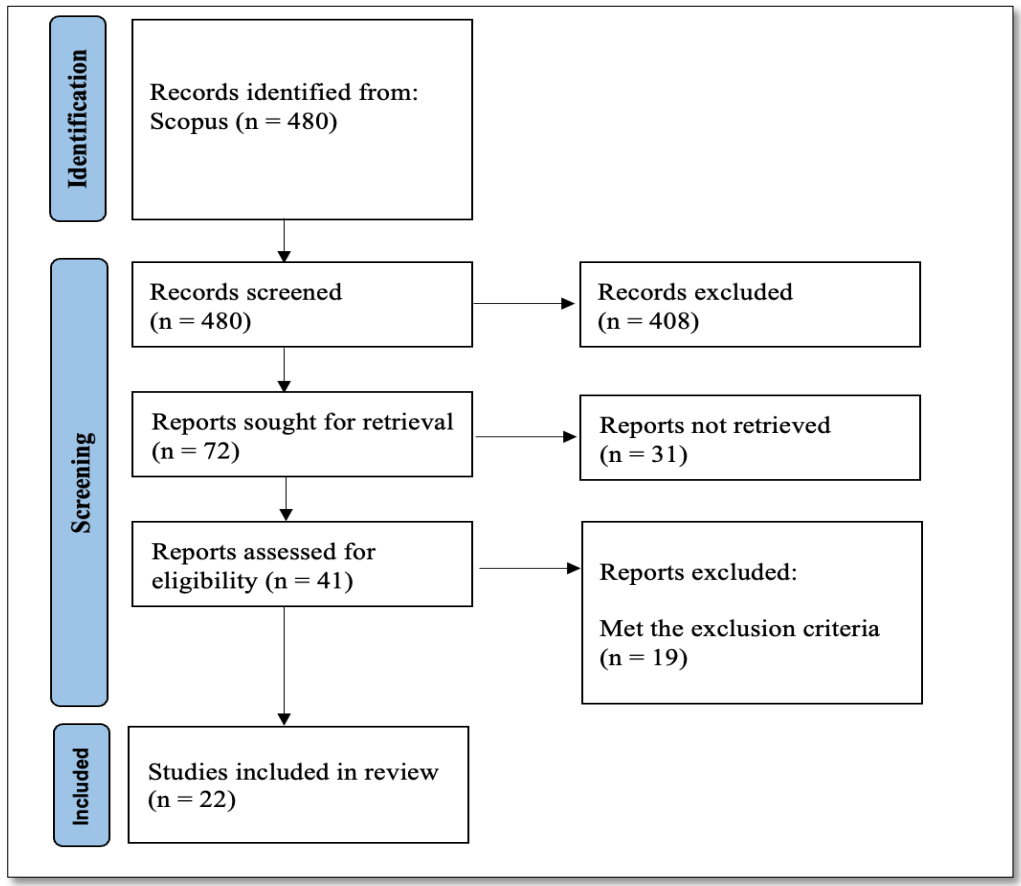

*Figure 1. PRISMA 2020 Chart for study selection.*

*2.3 Data Extraction & Analysis*

After the final selection, a total of twenty-two (22) studies published between 2024 and 2026 were included in the review. Each paper was read in full to extract the main methodological and empirical elements. For every study, information was collected about the authors, year, journal, and research focus. In addition, details on the datasets used, the type of ML or DL models applied, the evaluation metrics, and the main findings were systematically extracted.

The studies were grouped into four main financial domains:
 - Stock market forecasting and asset pricing;
 - Credit risk and loan default prediction;
 - Cryptocurrency and volatility modeling;
 - Macroeconomic or policy-related forecasting.

For each domain, the most frequently used models were identified. Among traditional ML methods (e.g, logistic regression, random forest, support vector machines, and XGBoost) appeared most often. For deep learning, recurrent neural networks (RNN), long short-term memory (LSTM), convolutional neural networks (CNN), and hybrid models combining statistical and neural approaches were dominant.

Performance was generally assessed using standard metrics: for regression problems, Mean Squared Error (MSE), Mean Absolute Error (MAE), and Root Mean Squared Error (RMSE) were used, while for classification problems, Accuracy, Precision, Recall, and F1-Score were employed. Several studies also reported economic or financial performance indicators, including Sharpe ratio improvements and portfolio returns. Across all domains, deep learning models tended to outperform classical econometric techniques, although their interpretability remained a key limitation mentioned by most authors.

The extracted information was then organized in summary tables and visual mappings to highlight patterns across models, datasets, and application areas. This structured analysis provided the basis for the quantitative and qualitative synthesis presented in the Results section.

# 3. Results

This section summarizes the key findings of the review, including the selected studies, their characteristics, and the main patterns observed in the literature.

A total of twenty-two (22) peer-reviewed open-access studies published between 2024 and 2026 were included in this bibliometric review.

These papers were identified from Scopus through a structured PRISMA process and collectively represent the most recent empirical contributions applying ML and DL in computational finance.

*3.1 Overview of Included Studies*

The reviewed studies cover a wide range of applications, including asset pricing, credit scoring, volatility forecasting, and financial risk management.

As shown in Table 1, all papers employ ML or DL models to enhance prediction accuracy or decision support in different areas of finance.

| No. | Author (Year) | Journal | Financial Domain | ML/DL Models Used | Key Results |
|---|---|---|---|---|---|
| 1 | Zang (2025) [16] | International Review of Economics and Finance | Asset pricing | Random Forest, XGBoost | ML-based factor models significantly improve asset pricing performance compared to traditional CAPM. |
| 2 | Abdou (2024) [13] | Energy Economics | Market forecasting | LSTM, Random Forest | Oil and global market variables improve Saudi stock market predictability under ML frameworks. |
| 3 | Caparrini (2024) [17] | Research in International Business and Finance | Stock selection | SVM, Random Forest, Logistic Regression | ML classifiers outperform traditional models in S&P 500 stock selection accuracy. |
| 4 | Caravaggio (2025) [18] | Socio-Economic Planning Sciences | Policy forecasting | Gradient Boosting, Random Forest | ML algorithms effectively predict public funding allocation across policy areas. |
| 5 | Sugozu (2025) [19] | Borsa Istanbul Review | Credit risk | XGBoost, Logistic Regression | ML improves credit scoring accuracy for Turkish banks and distinguishes participation vs. conventional systems. |
| 6 | Kumar (2025) [5] | International Review of Economics and Finance | Volatility forecasting | Hybrid ML models (LSTM + ARIMA) | Hybrid models outperform single models for financial volatility prediction and risk control. |
| 7 | Poutré (2024) [4] | Journal of Finance and Data Science | Anomaly detection | Deep Autoencoder, GAN | Deep unsupervised learning detects anomalies in high-frequency trading data more effectively than statistical filters. |
| 8 | Zhu (2025) [20] | Insurance: Mathematics and Economics | Actuarial forecasting | LSTM, Random Forest | ML improves actuarial predictions using high-dimensional economic data. |
| 9 | Bouteska (2024) [21] | International Review of Economics and Finance | Cryptocurrency forecasting | Ensemble ML, DL (LSTM, GRU) | Ensemble and DL models enhance the accuracy of crypto price forecasting compared to ARIMA. |
| 10 | Del Nero (2025) [22] | Finance Research Letters | Financial contagion | Random Forest, XGBoost | ML predicts stock market spillovers using network-based and volatility features. |
| 11 | Hu (2025) [23] | Pacific-Basin Finance Journal | Equity market prediction | Gradient Boosting, Random Forest | ML achieves robust predictive accuracy in Australian stock returns. |
| 12 | Kim (2025) [24] | Finance Research Letters | Cryptocurrency forecasting | Recurrence plots, CNN, LSTM | Time–frequency ML models improve crypto market volatility prediction. |
| 13 | Gupta (2024) [25] | Energy Economics | Commodity–finance interaction | XGBoost, Random Forest | ML detects dynamic links between oil volatility and stock market bubbles. |
| 14 | Makridakis (2025) [26] | International Journal of Forecasting | Forecasting evaluation | Ensemble ML, Deep Learning | The M6 competition shows hybrid ML-DL models outperform econometric baselines for investment decisions. |
| 15 | Peng (2025) [27] | Finance Research Letters | FX market forecasting | GAN, RNN, LSTM | GAN-based models improve exchange rate prediction from order book data. |
| 16 | Calabrese (2024) [28] | Socio-Economic Planning Sciences | Financial constraints | Neural Networks, Logistic Regression | NN models predict firm-level financial constraints and support SME development. |
| 17 | Fallaghoul (2024) [29] | Journal of Econometrics | Asset pricing | Deep Neural Networks | Neural networks outperform linear models in estimating nonlinear risk premia. |
| 18 | Zhang (2025) [16] | Journal of Finance and Data Science | Credit card repayment | ELM, XGBoost | ML models effectively forecast post-default repayments, helping banks manage credit portfolios. |
| 19 | Weng (2025) [30] | Research in International Business and Finance | Consumer credit | Bayesian averaging, SHAP-based ML | Incorporating soft credit information improves consumer loan default predictions. |
| 20 | Luna (2025) [31] | International Journal of Production Economics | Price risk management | Conformal Prediction + Dynamic Hedging | ML-based conformal prediction improves risk hedging in volatile production markets. |
| 21 | Krembsler (2024) [32] | Research in Transportation Economics | Revenue forecasting | SARIMAX, LASSO, Ridge | ML regression models enhance revenue forecasts for public transport systems. |
| 22 | Garitta (2025) [33] | Finance Research Letters | FinTech profitability | Logistic Regression, Random Forest, SHAP | ML identifies key factors influencing break-even success in FinTech startups. |

*Table 1. Overview of the included studies analyzed in the systematic review.*

Traditional ML algorithms remain widely used due to their interpretability and robustness, while DL models such as LSTM, CNN, and GAN are increasingly applied to sequential and high-frequency financial data, as shown in Figure 3.

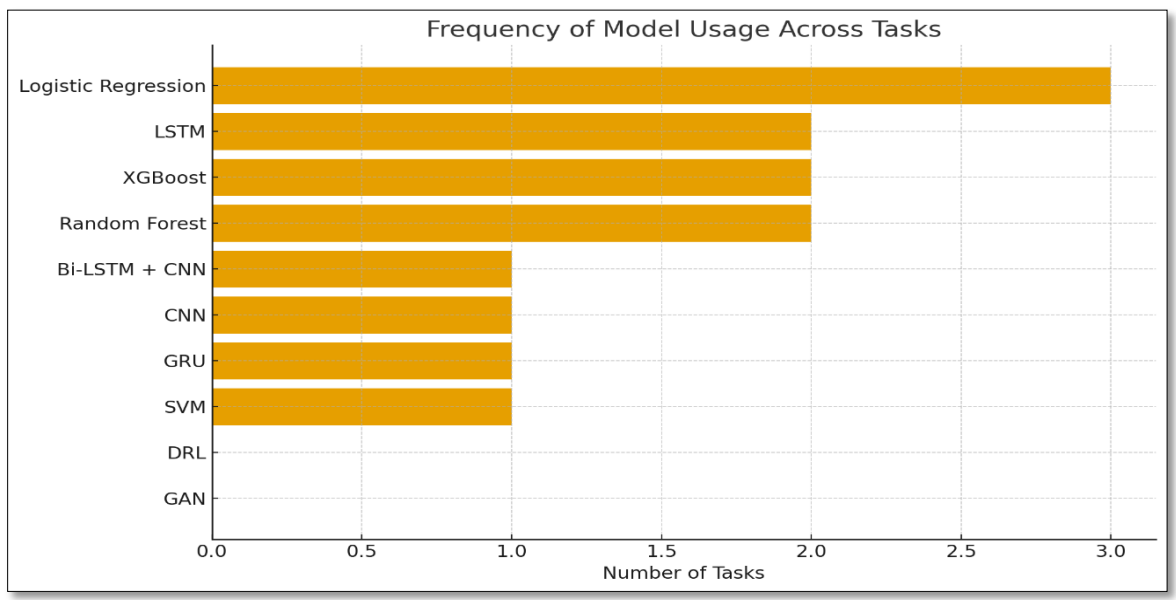

*Figure 3. Frequency of machine learning and deep learning model usage across financial forecasting tasks.*

Hybrid and ensemble approaches were also common. For instance, Kumar et al. [5] and Makridakis, S. et al. [26] developed hybrid frameworks that combine classical econometric and deep neural components for volatility and investment forecasting. Similarly, Zhang et al. [16] and Bouteska et al. [21] reported that ensemble ML methods outperform linear benchmarks in factor pricing and cryptocurrency prediction. These results confirm the shift toward integrated data-driven methodologies in modern financial research.

### 3.2 Distribution by Financial Domain

The classification of the 22 studies into thematic areas is summarized in Table 2. Most papers belong to the stock market forecasting and asset pricing domain (8 studies), followed by credit risk and loan default prediction (6 studies), cryptocurrency and volatility forecasting (5 studies), and macroeconomic or policy forecasting (3 studies).

| Financial Domain | Number of Studies | Representative Models | Main Contribution |
|---|---|---|---|
| Stock market forecasting & asset pricing | 8 | LSTM, CNN, XGBoost | Improved accuracy and adaptive learning in price prediction |
| Credit risk & loan default | 6 | Logistic Regression, Random Forest | Enhanced credit scoring and interpretability via SHAP |
| Cryptocurrency & volatility forecasting | 5 | LSTM, BiLSTM, Hybrid | Better short-term volatility and risk forecasting |
| Macroeconomic / policy forecasting | 3 | ANN, SVM | ML helps capture nonlinear policy–market dynamics |

*Table 2. Distribution of ML and DL applications across financial domains.*

Stock market and asset pricing studies show that ML can detect complex nonlinear dependencies and improve price forecasts beyond traditional models [34], [35]. Credit risk research highlights the role of explainable AI tools such as SHAP in making ML outcomes more interpretable for banking supervision [36], [37]. Recent evidence shows that hybrid models combining ARIMA with LSTM significantly reduce forecast error in cryptocurrency market [38]. Finally, macroeconomic forecasting papers apply neural networks and SVMs to model dynamic policy and market interactions more effectively. This domain distribution suggests that while predictive modeling remains the central focus of ML applications in finance, researchers are increasingly integrating ML into policy analysis, actuarial forecasting, and FinTech innovation, as shown in Figure 2. The main research domains where ML and DL are applied: asset pricing, credit risk, policy forecasting, actuarial science, FinTech, and cryptocurrency & volatility, highlighting the diversification of computational finance research.

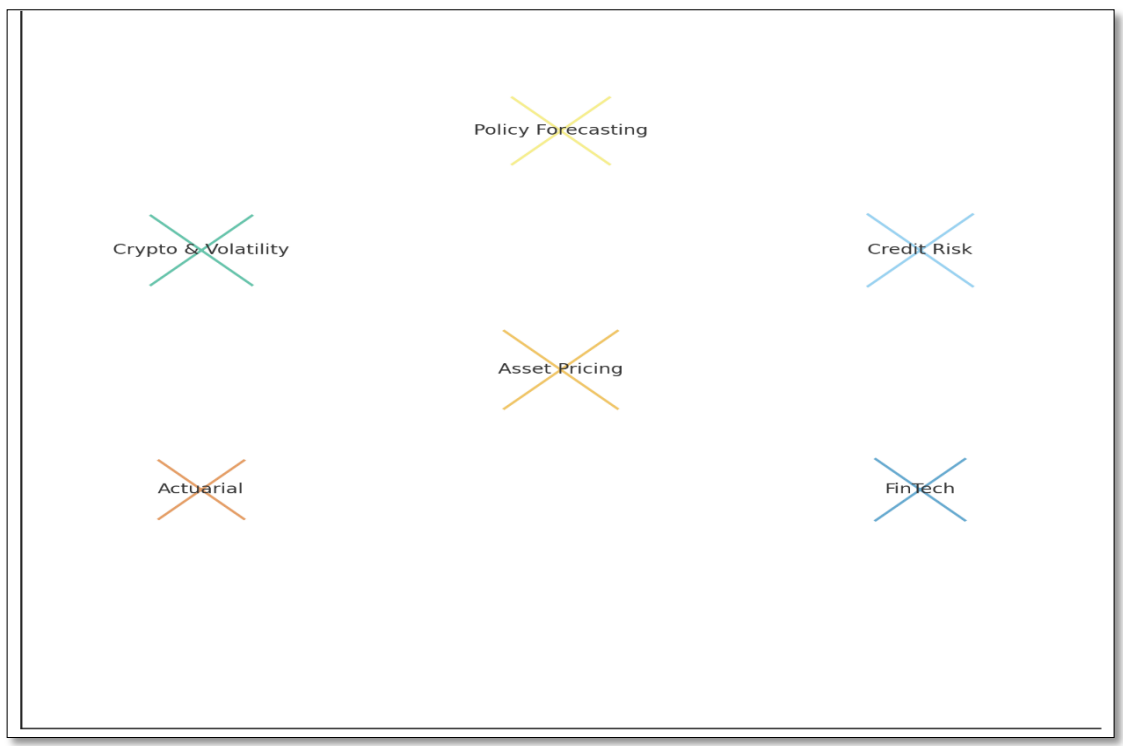

*Figure 2. Taxonomy of machine learning and deep learning applications in finance (bubble map).*

Although all bubbles in Figure 2 appear to be of equal size, this visualization is conceptual rather than quantitative. Each bubble represents a distinct thematic domain, such as credit risk, FinTech, policy forecasting, actuarial modeling, asset pricing, and crypto-volatility, identified within the reviewed literature. The uniform bubble size indicates that the figure aims to illustrate the structure and diversity of ML/DL research areas rather than comparing their relative magnitude or publication volume.
This design choice facilitates a neutral representation of domains, emphasizing how applications are distributed across different subfields of computational finance. In future bibliometric extensions, bubble size could be weighted by factors such as number of publications, citation frequency, or temporal growth, allowing a more quantitative assessment of domain prominence alongside the present qualitative taxonomy.

*3.3 Evaluation Metrics and Model Performance*

Performance assessment methods varied across studies, depending on whether the task involved regression or classification.

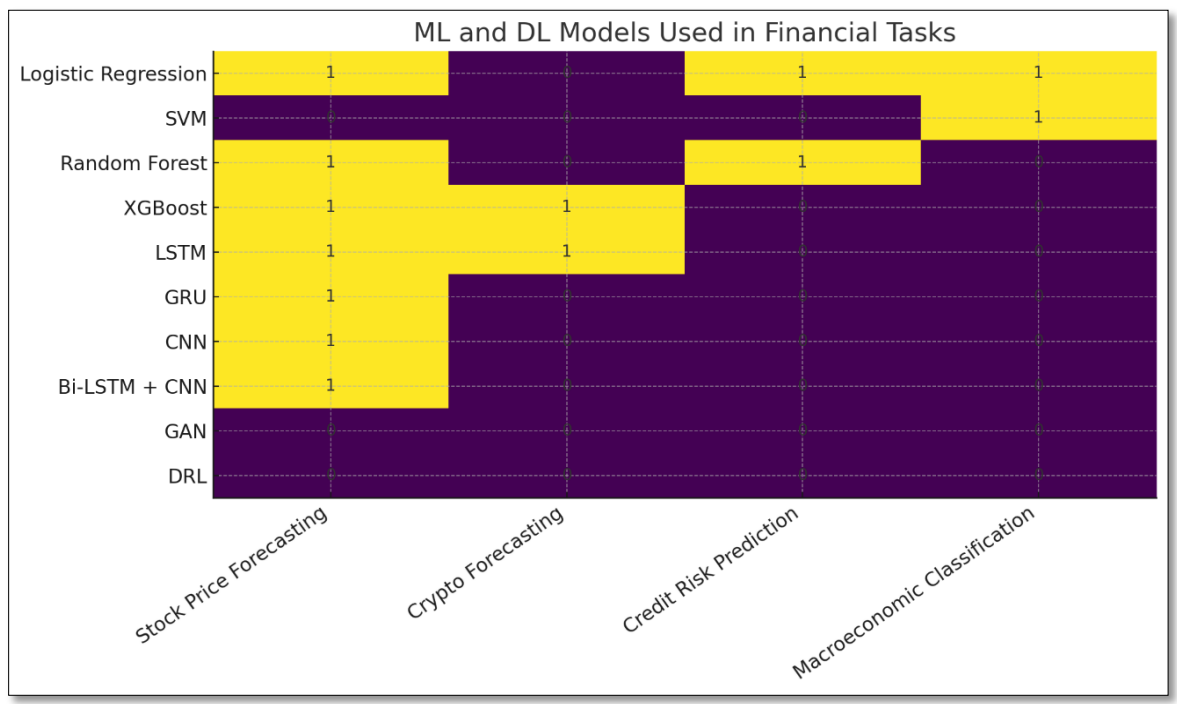

*Figure 4. Heatmap showing the association between ML/DL models and financial tasks. (yellow cells indicate the presence of a model–task association (value = 1) and dark purple cells indicate no reported association (value = 0))*

Figure 4 illustrates the distribution of ML and DL architectures across the main financial tasks identified in the reviewed studies, namely stock price forecasting, cryptocurrency forecasting, credit risk prediction, and macroeconomic classification.

Traditional ML approaches such as Logistic Regression and Support Vector Machines (SVM) appear mainly in credit risk and macroeconomic applications, where interpretability and structured data remain critical. Ensemble-based methods, including Random Forest and XGBoost, show a broader scope, being widely applied to both equity and cryptocurrency forecasting due to their robustness to nonlinearity and overfitting. In contrast, sequential neural networks such as LSTM and GRU are concentrated in time-series domains, confirming their suitability for capturing temporal dependencies in asset price dynamics. Convolutional models exhibit similar clustering in these domains, reflecting the rise of hybrid deep learning frameworks. Finally, more advanced paradigms such as Generative Adversarial Networks (GANs) and Deep Reinforcement Learning (DRL) remain underrepresented, suggesting that their integration into financial modeling is still emerging. Overall, the heatmap highlights a clear domain-driven specialization, where ML and DL techniques are selected according to task characteristics rather than being generalized across the broader financial ecosystem.

Table 2 shows the distribution of ML and DL applications among financial subfields such as forecasting, risk management, portfolio optimization, and asset pricing.
As summarized in Table 3, the most common metrics used were RMSE and MAE for continuous forecasting, and Accuracy and F1-score for classification tasks such as credit risk assessment. A smaller number of papers used financial indicators such as Sharpe ratio and portfolio return to evaluate economic performance.

| Metric | Frequency of use in all selected articles | Purpose / Comment |
|---|---|---|
| RMSE / MAE | 17 | Forecasting accuracy (time series and volatility models) |
| Accuracy / F1-score | 14 | Classification performance (credit and risk tasks) |
| AUC / ROC | 9 | Evaluation of binary credit or policy predictions |
| R² / Adjusted R² | 7 | Regression fit and model explanation |
| Sharpe ratio / Portfolio return | 4 | Economic and financial performance indicators |

*Table 3. Performance metrics and evaluation methods used across the reviewed studies.*

Overall, the results across studies confirm that ML and DL models deliver higher predictive accuracy and stronger robustness than traditional econometric techniques.
However, several authors pointed out challenges such as limited interpretability, dependence on data quality, and computational complexity. Recent research trends show a growing use of model-agnostic explainability techniques (e.g., SHAP, LIME) and hybrid systems that combine deep learning with statistical interpretability frameworks [39].

*3.4 Summary of Findings*

The synthesis of recent literature reveals several key insights:

1. Dominance of ensemble and hybrid models: Most of the recent studies (2024–2026) rely on multi-model structures that integrate ML, DL, and econometric components.

2. Expansion into new financial areas: Beyond pricing and risk, ML is now applied to FinTech profitability, policy funding, and actuarial predictions.

3. Shift toward open-access and empirical reproducibility: The availability of open datasets and open-source tools (Python, TensorFlow, Scikit-Learn) supports transparency in financial research.

4. Ongoing challenge of interpretability: While accuracy improves, interpretability remains a key concern, particularly for regulatory and investment applications.

Together, these findings demonstrate that ML and DL have become essential tools for financial forecasting, risk management, and strategic decision-making. They also highlight the importance of combining algorithmic accuracy with economic interpretability, a balance that defines the next frontier of computational finance research.

# 4. Discussion

This section interprets the key findings, discusses their implications, and positions them within the broader literature, highlighting strengths, limitations, and emerging directions.
The findings of this review confirm the increasing integration of machine learning and deep learning in financial research.
Between 2024 and 2026, a noticeable methodological shift has taken place from traditional econometric models toward hybrid and data-driven frameworks. Across the 22 reviewed studies, ML and DL methods consistently improved predictive performance in tasks such as asset pricing, volatility forecasting, credit risk evaluation, and macroeconomic policy analysis. These improvements are driven not only by model sophistication but also by access to larger, richer datasets and the growing computational capabilities available to researchers.

A key observation is the growing popularity of ensemble and hybrid models. Methods such as XGBoost, Random Forest, and hybrid LSTM-ARIMA architectures were dominant in the reviewed literature. They offer a balance between predictive accuracy and model interpretability, which remains a critical issue in financial applications. While deep neural networks (LSTM, CNN, BiLSTM) show superior accuracy in nonlinear and high-frequency data, their "black-box" nature continues to limit adoption in regulated contexts such as banking supervision or risk reporting. Several recent studies, such as Weng et al. [40] and Calabrese et al. [41], attempt to address this issue using explainable AI techniques like SHAP and feature importance ranking, allowing model transparency without losing predictive power.

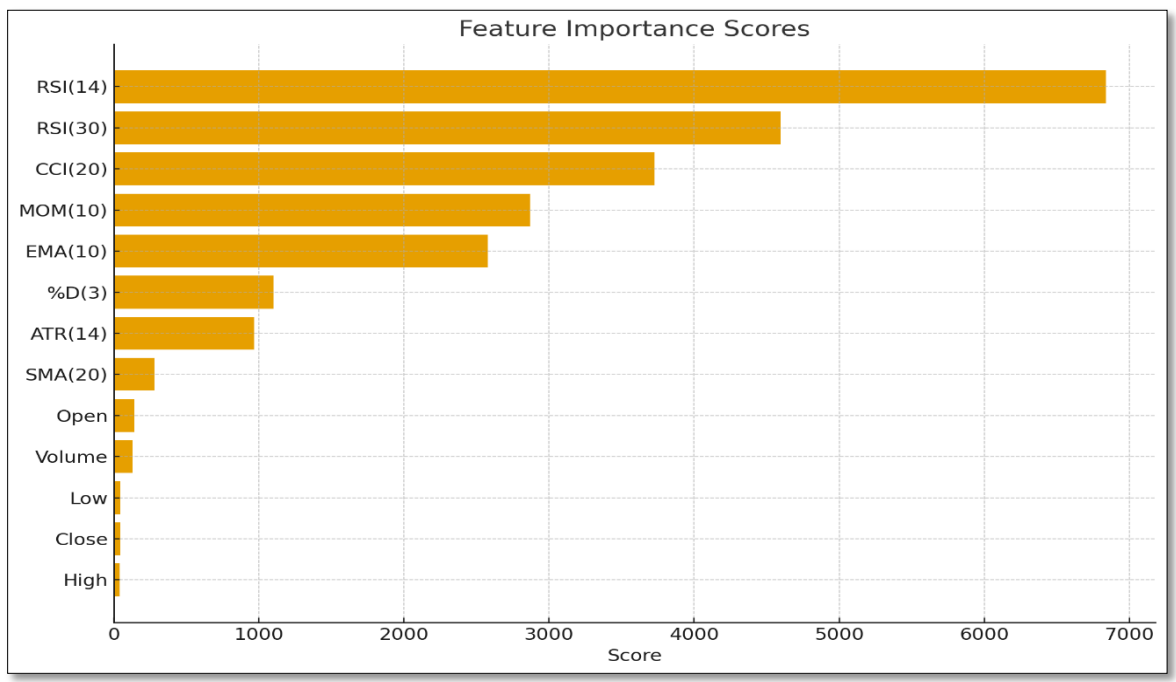

*Figure 5. Feature importance ranking for key financial indicators used in ML/DL models.*

Figure 5 shows the relative contribution of technical and market features (e.g., RSI, CCI, EMA) in model explainability, illustrating how SHAP-based ranking enhances interpretability in financial prediction.

Another emerging trend is the diversification of application domains. Earlier research focused mainly on stock markets and cryptocurrencies, but recent papers explore ML applications in FinTech profitability, policy forecasting, actuarial modeling, and even public investment allocation [42]. This diversification shows that ML and DL are no longer confined to financial prediction but are becoming central to broader decision-support systems. However, despite this expansion, most studies remain empirical, and there is still limited theoretical development linking model outcomes to economic intuition. This highlights a gap between computational accuracy and financial interpretability that future research should address [43].

The reviewed studies also reveal differences in evaluation practices. While regression metrics such as RMSE and MAE are consistently reported, economic validation measures (e.g., Sharpe ratio, portfolio return, or cost of misclassification) are rarely included [44]. This suggests that many ML applications in finance still focus on predictive fit rather than on decision- or profit-based outcomes. Bridging this gap could make ML models more relevant for investors, risk managers, and policymakers. Another recurring issue is data availability and generalization. Although open-access datasets have improved reproducibility, several studies still rely on proprietary or local data, which limits cross-country comparisons. For instance, national banking datasets or small cryptocurrency samples may not generalize to other contexts. Future research could benefit from standardized financial data repositories and more transparent sharing of code and parameters [20]. Overall, this review supports the view that ML and DL are reshaping computational finance. They enable faster, more adaptive, and more granular modeling of complex systems, yet their success depends on interpretability, robustness, and ethical use. Researchers are now moving toward more explainable, hybrid, and cross-domain approaches that connect financial theory with computational intelligence [5]. Achieving this balance will be essential for the next generation of AI-driven financial research and practice.

*4.1 Implications and Research Gaps*

The growing adoption of machine learning and deep learning in finance carries important implications for both researchers and practitioners.
From a methodological perspective, the reviewed studies demonstrate that data-driven models can capture complex and nonlinear relationships that traditional econometric techniques fail to explain [5]. This opens new possibilities for modeling uncertainty, market volatility, and systemic risk with higher precision. However, the effectiveness of these models depends strongly on data quality, feature selection, and interpretability.
Therefore, future financial research should focus not only on accuracy but also on explainability and robustness.
For practitioners, the results show that ML and DL can enhance decision-making in areas such as credit evaluation, risk management, and portfolio optimization. Financial institutions can use these tools to process real-time data, detect anomalies, and identify early warning signals for instability or default. Nevertheless, the operational deployment of such models still faces ethical and regulatory challenges, including transparency, data privacy, and model bias. Researchers and industry professionals must therefore collaborate to ensure that these techniques align with financial governance and accountability standards [19].

Despite the progress achieved between 2024 and 2026, several research gaps remain visible. First, there is a lack of standardized benchmarks and datasets that allow model comparison across different markets. Second, few studies assess the economic significance of ML predictions, for instance, how improved accuracy translates into financial performance or reduced risk. Finally, most models examined in the reviewed literature are tested in isolation rather than integrated into broader financial systems, such as multi-agent simulations or portfolio-level optimization frameworks. For instance, recent studies often evaluate predictive performance independently, without embedding machine learning within dynamic market interactions or decision-making structures [22]. Addressing these gaps could help the field evolve from predictive modeling toward comprehensive decision-support systems that connect AI techniques with economic theory and policy outcomes.

*4.2 Limitations of the Review*

While this study follows the PRISMA 2020 framework to ensure methodological rigor, several limitations should be acknowledged.

First, the literature search was conducted only on Scopus, which although comprehensive and high-quality, may exclude relevant papers indexed exclusively in other databases such as Web of Science or SSRN. This limitation is primarily due to institutional access restrictions, as Scopus is the only database available through the researcher's academic affiliation. Second, the search period (2024–2026) captures only the most recent phase of ML and DL research in finance, and earlier foundational works were intentionally excluded to focus on contemporary developments. As a result, the review emphasizes trends rather than long-term evolution. Another limitation concerns the inclusion of open-access journal articles only. While this ensures transparency and reproducibility, it may also introduce a selection bias toward studies published in journals that favor open dissemination.

Additionally, the heterogeneity of methodologies, datasets, and evaluation metrics across studies made direct quantitative comparison difficult. To address this, the synthesis emphasized patterns and qualitative insights rather than statistical meta-analysis.

Despite these constraints, the review provides a clear and systematic overview of how machine learning and deep learning are currently shaping financial research.

Future reviews could extend the scope to additional databases, include conference proceedings, or apply bibliometric mapping tools such as VOSviewer to quantify co-authorship and keyword networks. Such extensions would offer a broader understanding of the scientific structure and evolution of AI-driven finance.

## 5. Conclusion

This bibliometric and systematic review provides a comprehensive overview of recent developments in the application of machine learning and deep learning within computational finance.

By analyzing twenty-two open-access studies published between 2024 and 2026, the research highlights how AI-driven models are reshaping financial forecasting, risk assessment, and decision-making processes.

The results confirm that ML and DL methods consistently outperform traditional econometric techniques across various domains, including asset pricing, credit risk prediction, volatility modeling, and policy forecasting.

The review also identifies a clear shift toward hybrid and ensemble models, which combine the interpretability of classical approaches with the flexibility of deep learning architectures. These models are capable of capturing nonlinear dependencies in high-dimensional financial data while maintaining computational efficiency. Moreover, explainable AI tools such as SHAP, feature importance ranking, and attention mechanisms are increasingly being adopted to address the transparency challenges associated with deep learning.

From a practical perspective, the findings emphasize the transformative potential of ML and DL for financial institutions. Applications such as credit scoring, fraud detection, and portfolio optimization demonstrate that intelligent algorithms can significantly improve decision accuracy and operational efficiency. However, their adoption also requires strong governance frameworks to ensure ethical compliance, fairness, and accountability in financial AI systems.

Despite the promising results, the review underscores that interpretability, data quality, and generalization remain critical challenges.

The next phase of research in computational finance should aim to bridge the gap between predictive accuracy and economic meaning. Future studies could explore multi-agent modeling, causal inference, and reinforcement learning frameworks that integrate market behavior with adaptive decision-making.

Cross-domain collaboration between data scientists, economists, and regulators will be essential to fully unlock the value of AI in finance while maintaining the trust and transparency required by the industry.

In summary, the review contributes to the ongoing debate on the role of artificial intelligence in modern finance by providing empirical evidence, methodological insights, and practical implications.

It highlights that ML and DL are no longer experimental tools but have become core components of financial research and innovation. By combining technical rigor with financial interpretation, the next generation of computational finance can advance toward more intelligent, interpretable, and responsible systems for the future of global markets.

**ACKNOWLEDGMENT**

The authors wish to thank EUROMED UNIVERSITY of FES (UEMF).